\documentclass[11pt]{article}
\usepackage{amsmath, amssymb, graphicx}
\newcommand{\qed}{\hfill\vrule height6pt  width6pt
depth0pt \medskip}
\title{Shortest paths in the Tower of Hanoi graph and finite automata}
\author{Dan Romik
\thanks{Department of Mathematics, The Weizmann Institute of Science,
Rehovot 76100, Israel. email: \tt{romik@wisdom.weizmann.ac.il}}
}
\begin{document}
\maketitle
\begin{abstract}
We present efficient algorithms for constructing a shortest path
between two states in the Tower of Hanoi graph, and for computing the
length of the shortest path. The key element is a finite-state machine
which decides, after examining on the average only $\frac{63}{38}
\approx 1.66$ of the largest discs, whether the largest disc will be
moved once or twice. This solves a problem raised by Andreas Hinz,
and results in a better understanding of how the shortest path is
determined. Our algorithm for computing the length of the shortest
path is typically about twice as fast as the existing algorithm. We
also use our results to give a new derivation of the average distance
$\frac{466}{885}$ between two random points on the Sierpi\'nski gasket
of unit side.
\end{abstract}

\subsection*{1. Introduction}

The \emph{Tower of Hanoi} puzzle, invented in 1883 by the French
mathematician Edouard Lucas, has become a classic example in the
analysis of algorithms and discrete mathematical structures (see
e.g. \cite{concmath}, 1.1). The puzzle consists of $n$ discs of
different sizes, stacked on three vertical pegs, in such a way that no
disc lies on top of a smaller disc. A permissible \emph{move} is to
take the top disc from one of the pegs and move it to one of the other
pegs, as long as it is not placed on top of a smaller disc. The set of
states of the puzzle, together with the permissible moves, thus forms
a graph in a natural way. The number of vertices in the $n$-disc Hanoi
graph is $3^n$.

The main question of interest is to find \emph{shortest paths} in the
state graph, i.e., shortest sequences of moves leading from a given
initial state to a given terminal state. The simplest and most
well-known case is that in which it is required to move all the discs
from one of the pegs to another, i.e. where the initial and terminal
states are two of the three ``perfect'' states with all the discs on
the same peg. This is very easy, and can be shown to take exactly
$2^n-1$ moves. More difficult is to get from a given arbitrary initial
state to one of the perfect states - Hinz \cite{hinz1} calls this the
``p1'' problem. This takes $2^n-1$ moves in the worst case (which is
when the initial state is another perfect state), and on the average
$\frac{2}{3}\cdot (2^n-1)$ moves for a randomly chosen initial state
\cite{er}. Moreover, there is a simple and efficient algorithm to
compute the shortest path in this case.

In the most general case of arbitrary initial \emph{and} terminal
states, however, the question of computing the shortest path and its
length (the ``p2'' problem \cite{hinz1}) in the most efficient manner,
has not been completely resolved so far. (The worst-case behavior is
still $2^n-1$ moves, and the average number of moves for random
initial and terminal states has been shown \cite{chan},\cite{hinz3} to
be about $\frac{466}{885}\cdot 2^n$.) The main obstacle in the
understanding of the behavior of the shortest path, has been the
behavior of the largest disc that ``separates'' the initial and
terminal states, i.e. the largest disc which is not on the same peg in
both states (trivially, any larger discs may simply be ignored). It is
not difficult to see \cite{hinz1} that in a shortest path, this disc
will be moved either once (from the source peg to the target peg) or
twice (from the source to the target, via the third peg). The problem
is to decide which of the two alternatives is the correct one. Once
this is settled, the path may be constructed by two applications of
the algorithm for the p1 problem. Hinz \cite{hinz1} proposed an
algorithm for the computation of the shortest path based on this
idea. The algorithm consists essentially of computing the length of
the path for both alternatives and choosing the shorter of the two.

In this paper, we propose a more thorough explanation of the process
whereby it is decided which of the two paths is the shortest. We show
that it is possible to keep track of the relevant information using a
finite-state machine, which at each step reads the locations of the
next-smaller disc in the initial and terminal states, and changes its
internal state accordingly. Eventually, the machine reaches a terminal
state, whereupon it pronounces which of the two paths is the
shortest. For a random input, its expected stopping time is computed
to be $\frac{63}{38}$. In other words, after observing on the average
the locations of just the $1.66$ largest discs in the initial and
terminal states, we will know which of the paths to choose, and we
will be able to continue using the algorithm for the p1 problem. If
one is interested just in the length of the shortest path, then our
algorithm is typically about twice as fast as the algorithm proposed
by Hinz \cite{hinz1} (with a small constant overhead due to the
initial 1.66 discs), since it overrides the need to compute both the
distance for the path that moves the largest disc once, and the path
that moves it twice. 

The paper is organized as follows: In the next section, we define the
\emph{discrete Sierpi\'nski gasket} graph, a graph which is isomorphic
to the Tower of Hanoi state graph, but for which the labeling of the
vertices is simpler to understand. In section 3, we present the main
ideas for this graph, and then in section 4 show how to translate the
results to the Hanoi graph by a re-labeling of the vertices. In
section 5 we perform a probabilistic analysis of the finite-state
machine, to compute the average number $\frac{63}{38}$ of discs that
need to be read in order to decide whether the largest disc will be
moved once or twice, and to give a new derivation of the asymptotic
value $\frac{466}{885}\cdot 2^n$ for the average distance between two
random states in the $n$-disc Hanoi graph.

\subsection*{2. The discrete Sierpi\'nski gasket}

We now define a family of graphs called \emph{discrete Sierpi\'nski
gaskets}. These graphs are finite versions of the famous fractal
constructed by the Polish mathematician Waclaw Sierpi\'nski in 1915. The
connection between the Tower of Hanoi problem and the Sierpi\'nski
gasket was first observed by Ian Stewart \cite{stewart}, and was later
used by Andreas Hinz and Andreas Schief \cite{hinzschief} in their
calculation of the average distance between points on the Sierpi\'nski
gasket.

The $n$th discrete Sierpi\'nski gasket graph, which we denote by $SG_n$,
consists of the vertex set $V(SG_n) = \{ T, L, R \}^n$ (the symbols
$T,L,R$ indicate ``top'', ``left'' and ``right'', respectively), with
the edges defined as follows: First, for each $x = a_{n-1} a_{n-2}
... a_1 a_0 \in V(SG_n)$ (for reasons that will become apparent below,
this will be our standard indexing of the coordinates of the vertices
of $SG_n$) we have edges connecting $x$ to
$$ a_{n-1} a_{n-2} ... a_1 \beta, \qquad \beta \in \{T,L,R\}\setminus
\{a_0\} $$
Second, define the \emph{tail} of $x=a_{n-1} a_{n-2}... a_0$ as the
suffix $a_k a_{k-1} ... a_1 a_0$ of $x$, where $k$ is maximal such
that $a_k = a_{k-1} = ... = a_0$. If $x$ has a tail of length $k+1<n$,
then $x$ is of the form $a_{n-1} a_{n-2} ... a_{k+2} \beta \alpha
\alpha ... \alpha$. Connect $x$ with an edge to the vertex
$$ a_{n-1} a_{n-2} ... a_{k+2} \alpha \beta \beta ... \beta $$

One possible embedding of $SG_n$ in the plane is illustrated in Figure
1 below. This embedding makes clear the meaning of the labeling of the
vertices: The first letter (the ``most significant digit'') signifies
whether the vertex is in the top, left or right triangles inside the
big triangle; the next letter locates the vertex within the top, left
or right thirds of that triangle, etc.

\hspace{-50.0 pt}
\includegraphics{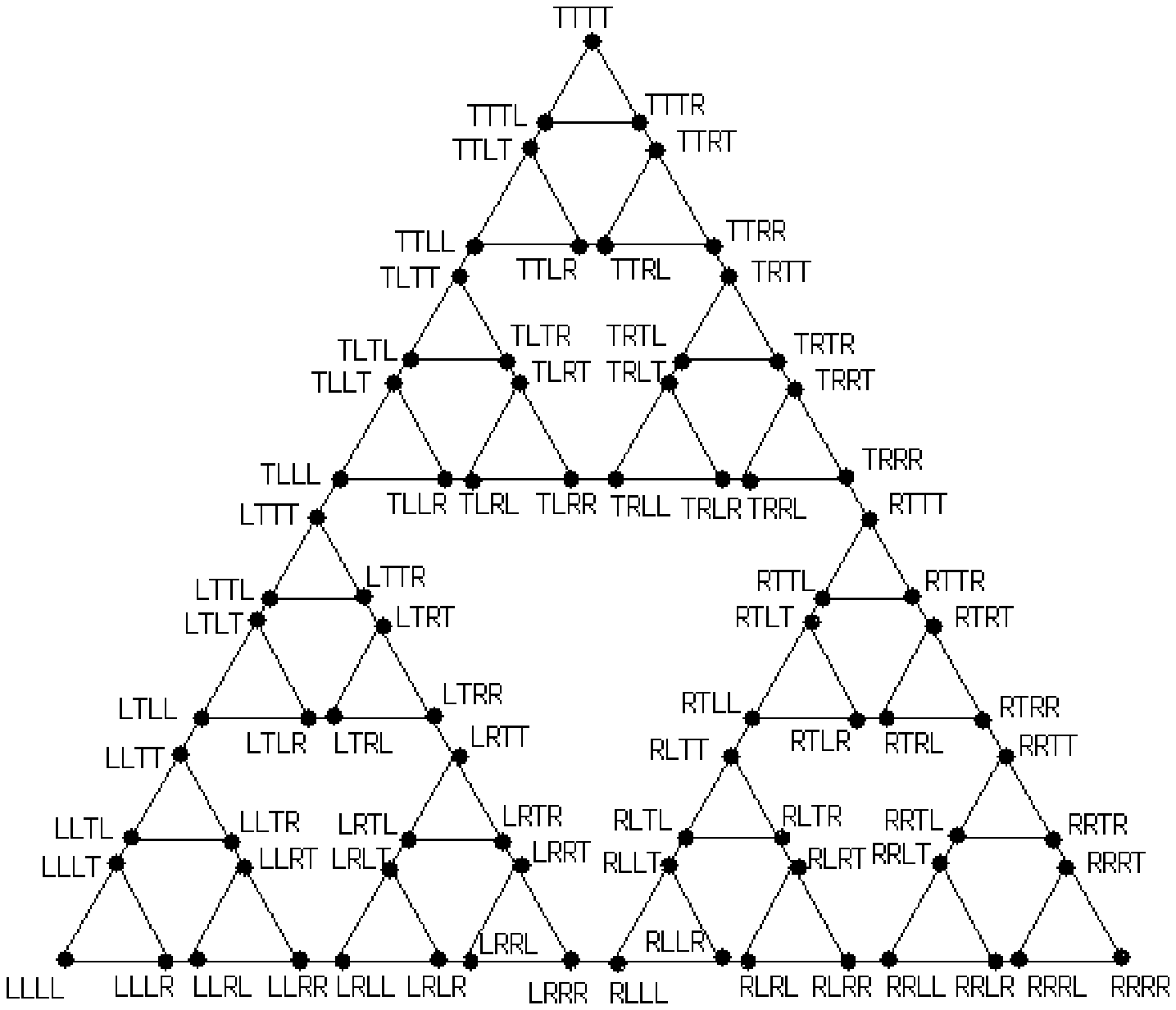}

\vspace{-20.0 pt}
\begin{center}
Figure 1: The graph $SG_4$
\end{center}

\bigskip
It will be shown in section 4 that $SG_n$ is isomorphic, in a
computationally straightforward way, to the $n$-disc Hanoi graph. (The
same was shown in \cite{hinzschief}, with less emphasis on explicit computation
of the isomorphism.) Thus, the problem of shortest paths on the Hanoi
graph reduces to that of shortest paths in the discrete Sierpi\'nski
gasket. We tackle this problem in the next section.

\subsection*{3. Shortest paths in $SG_n$}

For vertices $x,y \in V(SG_n)$, we define the distance $d(x,y)$ to be
the length of a shortest path from $x$ to $y$. Our goal is to write
down a recursion equation for this distance, which is at the heart of
the finite-state machine we will construct to compute $d(x,y)$. First,
let us review briefly some of the known facts about $d(x,y)$ in the
simple case when $y$ is one of the ``perfect'' states $LLL...L,
RR...R, TT...T$. For concreteness, assume that $y=LLL...L$, and let
$x=a_{n-1} a_{n-2}...a_1 a_0 \in V(SG_n)$ as before. Then it is known
that
$$ d(x,y) = \sum_{a_k \neq L} 2^k $$ 
A simple algorithm exists for computing a shortest path from $x$ to
$y$ in this case. In the Hanoi labeling of the graph, the algorithm is
described in \cite{hinz1}. In the current labeling, the algorithm is even
simpler and is based on the binary number system: if one identifies
the symbol $L$ with $0$ and the symbols $R$ and $T$ with $1$, then
traversing the edges of the graph becomes equivalent to the operations
of subtraction or addition of 1 in binary notation. The number of
steps to reach $LL...L \equiv 00...0$ is then clearly the right-hand
side in the above equation.

With these preparatory remarks, we now attack the problem of general
$x=a_{n-1}a_{n-2} ... a_0,\ \  y=b_{n-1}b_{n-2}...b_0$. First, observe
that we may assume that $a_{n-1}\neq b_{n-1}$, since otherwise we may
simply consider $x$ and $y$ as vertices in the graph $SG_{n-1}$ (note
the self-similar structure in the definition of the graph, also
apparent in the Tower of Hanoi puzzle when one ignores the largest disc).
For concreteness, we shall analyze in detail the case where
$a_{n-1}=T,\ \ b_{n-1}=R$. Referring to Figure 1 for convenience, we
see that
$$ d(x,y) = \min\bigg(1+d(x,TRRR...R)+d(y,RTT...T), $$ $$ 1+2^{n-1}+
d(x,TLLL...L) + d(y,RLL...L)\bigg), $$ 
since in a shortest path from $x$ to $y$, one must go from the top
triangle to the right triangle either through the edge $\{ TRR...R,
RTT...T \}$ (we call this Alternative 1, see Theorem 1 below) or
through a shortest path from $LTT...T$ to $LRR...R$ (Alternative 2) -
in the Tower of Hanoi language, this is an indication of the fact that
in a shortest sequence of moves the largest disc must move either once
or twice, see \cite{hinz1}.

To simplify the next few equations, introduce the following notation:
if $u=c_{n-1}c_{n-2}...c_0, \in \{T,L,R\}^n$, let $u' =
c_{n-2}c_{n-3}... c_0$, and define for any $\alpha \in \{ L, T, R \}$
$$ f_\alpha(u) = \sum_{ c_k \neq \alpha } 2^k $$
Then we have
$$ d(x,y) = 1 + \min\bigg( f_R(x')+f_T(y'),\ \  2^{n-1} + f_L(x')+f_L(y')
\bigg) $$
The recursion equations which will enable us to construct our
finite-state machine and compute $d(x,y)$ are now given by the
following theorem:

\paragraph{Theorem 1 - The finite state machine.} For
$u=c_{n-1}c_{n-2} ... c_0,\ \ v=d_{n-1}d_{n-2}...d_0 \in \{T,L,R\}^n$,
define the functions
$$ p(u,v) = \min\bigg(f_R(u)+f_T(v),\ 2^n + f_L(u)+f_L(v)\bigg) $$
$$ q(u,v) = \min\bigg(2^n + f_R(u)+f_T(v),\ f_L(u)+f_L(v)\bigg) $$
$$ r(u,v) = \min\bigg(f_R(u)+f_T(v),\ f_L(u)+f_L(v)\bigg) $$
(note that $p,q,r$ depend implicitly on the length $n$ of the strings.)
Then we have the equations
$$
p(u,v) = \left\{\begin{array}{lll}
  f_R(u)+f_T(v) & \begin{array}{l} c_{n-1}=R, d_{n-1}=T \textrm{  or}\\
                                   c_{n-1}=R, d_{n-1}=L \textrm{  or}\\
                                   c_{n-1}=R, d_{n-1}=R \textrm{  or}\\
                                   c_{n-1}=L, d_{n-1}=T \textrm{  or}\\
                                   c_{n-1}=T, d_{n-1}=T \textrm{  or}\\
                                   c_{n-1}=T, d_{n-1}=R
                  \end{array} & \textrm{(Alternative 1)}
     \\ \ \\
  2^n + p(u',v') & \begin{array}{l} 
                                   c_{n-1}=T,d_{n-1}=L \textrm{  or}\\
                                   c_{n-1}=L,d_{n-1}=R
                   \end{array} \\ \ \\
  2^n + r(u',v') & \begin{array}{l} c_{n-1}=L, d_{n-1}=L \end{array}
\end{array} \right.
$$
$$
q(u,v) = \left\{\begin{array}{lll}
  f_L(u)+f_L(v) & \begin{array}{l} c_{n-1}=L, d_{n-1}=L \textrm{  or}\\
                                   c_{n-1}=L, d_{n-1}=T \textrm{  or}\\
                                   c_{n-1}=L, d_{n-1}=R \textrm{  or}\\
                                   c_{n-1}=R, d_{n-1}=L \textrm{  or}\\
                                   c_{n-1}=T, d_{n-1}=L \textrm{  or}\\
                                   c_{n-1}=T, d_{n-1}=R
                  \end{array} & \textrm{(Alternative 2)}\\ \ \\
  2^n + q(u',v') & \begin{array}{l} 
                                   c_{n-1}=T,d_{n-1}=T \textrm{  or}\\
                                   c_{n-1}=R,d_{n-1}=R
                   \end{array} \\ \ \\
  2^n + r(u',v') & \begin{array}{l} c_{n-1}=R, d_{n-1}=T \end{array}
\end{array} \right.
$$
$$
r(u,v) = \left\{ \begin{array}{lll}
  f_R(u)+f_T(v) & \begin{array}{l} c_{n-1}=R, d_{n-1}=T \end{array} 
 & \textrm{(Alternative 1)} \\ \ \\
  f_L(u)+f_L(v) & \begin{array}{l} c_{n-1}=L, d_{n-1}=L \end{array} 
 & \textrm{(Alternative 2)} \\ \ \\
  2^{n-1}+r(u',v') & \begin{array}{l} c_{n-1}=L,d_{n-1}=T\textrm{  or}\\
                                      c_{n-1}=R,d_{n-1}=L
                     \end{array}\\
 \ \\
  2^n+r(u',v') & \begin{array}{l} c_{n-1}=T, d_{n-1}=R \end{array} \\
 \ \\
  2^{n-1}+p(u',v') & \begin{array}{l}
                      c_{n-1}=R, d_{n-1}=R \textrm{  or}\\
                      c_{n-1}=T, d_{n-1}=T
                     \end{array} \\
 \ \\
  2^{n-1}+q(u',v') & \begin{array}{l}
                      c_{n-1}=T, d_{n-1}=L \textrm{  or}\\
                      c_{n-1}=L, d_{n-1}=R
                     \end{array}
\end{array} \right.
$$
These equations will hold even for $n=1$ if one defines trivially
$p(u,v)=q(u,v)=r(u,v)=0$ for $u,v = \emptyset \in
\{T,L,R\}^0=\{\emptyset\}$. Alternatives 1,2 in the parentheses
signify whether the minimum is attained by its first or second
arguments, respectively.

\paragraph{Proof.} The proof consists simply of inspecting the
definitions of the functions $p,q,r$ and verifying that in each of the
cases the required relations hold. We omit the details, since the
reader would no doubt have to go through the same thought process to
verify them on her own as to check that a purported proof is correct.
\qed

\bigskip
A schematic representation of the finite state machine is shown in
Figures 2, 3 and 4 below. We present several variants of the machine:
the machine in Figure 2 only decides between Alternative 1 and
Alternative 2. The machine in Figure 3, which has auxiliary counters
for the distance and for the variable $n$ (so strictly speaking it is
not really a finite-state automaton), actually computes $d(x,y)$. Note
that this is still designed only for the case in which $x$ begins with
the symbol $T$ and $y$ begins with $R$. Figure 4 shows the complete
finite-state machine which may be used to compute $d(x,y)$ for any two
states $x,y \in V(SG_n)$. This includes an initial component that
discards the first few symbols which are identical for $x$ and $y$,
and another component that permutes the symbols $T,L,R$ to fit the
design of the basic machine in Figure 2.

\begin{picture}(300,230)(0,0)
\put(30,100){\circle{30}}
\put(150,100){\circle{30}}
\put(270,100){\circle{30}}
\put(75,190){\framebox(30,30)}
\put(195,190){\framebox(30,30)}
\put(50,90){\vector(1,0){80}}
\put(130,110){\vector(-1,0){80}}
\put(170,90){\vector(1,0){80}}
\put(250,110){\vector(-1,0){80}}
\qbezier(20,82)(0,55)(25,45)
\qbezier(40,82)(60,55)(35,45)
\qbezier(25,45)(30,43)(35,45)
\put(45,77){\vector(-1,1){7}}
\qbezier(140,82)(120,55)(145,45)
\qbezier(160,82)(180,55)(155,45)
\qbezier(145,45)(150,43)(155,45)
\put(165,77){\vector(-1,1){7}}
\qbezier(260,82)(240,55)(265,45)
\qbezier(280,82)(300,55)(275,45)
\qbezier(265,45)(270,43)(275,45)
\put(285,77){\vector(-1,1){7}}
\put(35,120){\vector(2,3){43}}
\put(145,120){\vector(-2,3){43}}
\put(155,120){\vector(2,3){43}}
\put(265,120){\vector(-2,3){43}}
\put(55,173){\shortstack{$RT$}}
\put(47,163){\shortstack{$RL$}}
\put(39,153){\shortstack{$RR$}}
\put(33,143){\shortstack{$LT$}}
\put(28,133){\shortstack{$TT$}}
\put(22,123){\shortstack{$TR$}}
\put(232,173){\shortstack{$LL$}}
\put(240,163){\shortstack{$LT$}}
\put(248,153){\shortstack{$LR$}}
\put(253,143){\shortstack{$RL$}}
\put(260,133){\shortstack{$TL$}}
\put(266,123){\shortstack{$TR$}}
\put(85,80){\shortstack{$LL$}}
\put(70,112){\shortstack{$RR,\ TT$}}
\put(195,80){\shortstack{$TL,\ LR$}}
\put(205,112){\shortstack{$RT$}}
\put(0,40){\shortstack{$TL,$}}
\put(25,32){\shortstack{$LR$}}
\put(240,40){\shortstack{$TT,$}}
\put(265,32){\shortstack{$RR$}}
\put(114,67){\shortstack{$LT$}}
\put(114,57){\shortstack{$RL$}}
\put(114,47){\shortstack{$TR$}}
\put(104,146){\shortstack{$RT$}}
\put(179,146){\shortstack{$LL$}}
\put(19,103){\shortstack{{\tiny START}}}
\put(18,94){\shortstack{\tiny (Alt. 1)}}
\put(141,94){\shortstack{\tiny (draw)}}
\put(258,94){\shortstack{\tiny (Alt. 2)}}
\put(81,208){\shortstack{{\tiny STOP}}}
\put(80,200){\shortstack{{\tiny Alt. 1}}}
\put(201,208){\shortstack{{\tiny STOP}}}
\put(200,200){\shortstack{{\tiny Alt. 2}}}
\end{picture}

\vspace{-20.0 pt}
\noindent
Figure 2: The finite state machine: deciding between Alternative 1 and
Alternative 2. The two letters signify the two inputs from $x$ and
$y$, reading at each step the next-most-significant symbol. The
parentheses in the non-terminal states indicate that if the input
terminates without a decision, then in the {\tiny START} state
Alternative 1 wins, in the rightmost state Alternative 2 wins, and in
the middle state there is a draw, meaning that the shortest path is
not unique and both alternatives are valid.

\newpage
\begin{picture}(300,200)(0,0)
\put(30,100){\circle{30}}
\put(150,100){\circle{30}}
\put(270,100){\circle{30}}
\put(75,190){\framebox(30,30)}
\put(195,190){\framebox(30,30)}
\put(50,90){\vector(1,0){80}}
\put(130,110){\vector(-1,0){80}}
\put(170,90){\vector(1,0){80}}
\put(250,110){\vector(-1,0){80}}
\qbezier(20,82)(0,55)(25,45)
\qbezier(40,82)(60,55)(35,45)
\qbezier(25,45)(30,43)(35,45)
\put(45,77){\vector(-1,1){7}}
\qbezier(137,85)(114,75)(105,60)
\qbezier(140,82)(130,59)(115,49)
\qbezier(105,60)(97,40)(115,49)
\put(137,77){\vector(2,3){5}}
\qbezier(163,85)(186,75)(195,60)
\qbezier(160,82)(170,59)(185,49)
\qbezier(195,60)(203,40)(185,49)
\put(164,77){\vector(-2,3){5}}
\qbezier(260,82)(240,55)(265,45)
\qbezier(280,82)(300,55)(275,45)
\qbezier(265,45)(270,43)(275,45)
\put(285,77){\vector(-1,1){7}}
\put(35,120){\vector(2,3){43}}
\put(145,120){\vector(-2,3){43}}
\put(155,120){\vector(2,3){43}}
\put(265,120){\vector(-2,3){43}}
\put(55,173){\shortstack{$RT$}}
\put(47,163){\shortstack{$RL$}}
\put(39,153){\shortstack{$RR$}}
\put(33,143){\shortstack{$LT$}}
\put(28,133){\shortstack{$TT$}}
\put(22,123){\shortstack{$TR$}}
\put(232,173){\shortstack{$LL$}}
\put(240,163){\shortstack{$LT$}}
\put(248,153){\shortstack{$LR$}}
\put(253,143){\shortstack{$RL$}}
\put(260,133){\shortstack{$TL$}}
\put(266,123){\shortstack{$TR$}}
\put(85,80){\shortstack{$LL$}}
\put(70,112){\shortstack{$RR,\ TT$}}
\put(195,80){\shortstack{$TL,\ LR$}}
\put(205,112){\shortstack{$RT$}}
\put(0,40){\shortstack{$TL,$}}
\put(25,32){\shortstack{$LR$}}
\put(240,40){\shortstack{$TT,$}}
\put(265,32){\shortstack{$RR$}}
\put(85,58){\shortstack{$LT$}}
\put(85,48){\shortstack{$RL$}}
\put(156,53){\shortstack{$TR$}}
\put(104,146){\shortstack{$RT$}}
\put(179,146){\shortstack{$LL$}}
\put(19,105){\shortstack{{\tiny START}}}
\put(21,98){\shortstack{\tiny $d=1$}}
\put(18,91){\shortstack{\tiny (Alt. 1)}}
\put(141,94){\shortstack{\tiny (draw)}}
\put(258,94){\shortstack{\tiny (Alt. 2)}}
\put(81,208){\shortstack{{\tiny STOP}}}
\put(80,200){\shortstack{{\tiny Alt. 1}}}
\put(201,208){\shortstack{{\tiny STOP}}}
\put(200,200){\shortstack{{\tiny Alt. 2}}}
\put(26,47){\shortstack{$2^n$}}
\put(106,37){\shortstack{$2^{n-1}$}}
\put(184,37){\shortstack{$2^n$}}
\put(93,91){\shortstack{$2^n$}}
\put(74,100){\shortstack{$2^{n-1}$}}
\put(213,91){\shortstack{$2^{n-1}$}}
\put(194,100){\shortstack{$2^n$}}
\put(266,47){\shortstack{$2^n$}}
\put(67,163){\shortstack{\tiny$f_R(x)+f_T(y)$}}
\put(187,163){\shortstack{\tiny$f_L(x)+f_L(y)$}}
\end{picture}

\vspace{-20.0 pt}
\noindent
Figure 3: The finite state machine: computing $d(x,y)$. Add to $d$ the
number on each edge traversed, decrease $n$ by 1 and replace $x$ by
$x'$ and $y$ by $y'$.

\hspace{30.0 pt}
\begin{picture}(300,340)(0,0)
\put(30,100){\circle{30}}
\put(150,100){\circle{30}}
\put(270,100){\circle{30}}
\put(75,190){\framebox(30,30)}
\put(195,190){\framebox(30,30)}
\put(50,90){\vector(1,0){80}}
\put(130,110){\vector(-1,0){80}}
\put(170,90){\vector(1,0){80}}
\put(250,110){\vector(-1,0){80}}
\qbezier(20,82)(0,55)(25,45)
\qbezier(40,82)(60,55)(35,45)
\qbezier(25,45)(30,43)(35,45)
\put(45,77){\vector(-1,1){7}}
\qbezier(137,85)(114,75)(105,60)
\qbezier(140,82)(130,59)(115,49)
\qbezier(105,60)(97,40)(115,49)
\put(137,77){\vector(2,3){5}}
\qbezier(163,85)(186,75)(195,60)
\qbezier(160,82)(170,59)(185,49)
\qbezier(195,60)(203,40)(185,49)
\put(164,77){\vector(-2,3){5}}
\qbezier(260,82)(240,55)(265,45)
\qbezier(280,82)(300,55)(275,45)
\qbezier(265,45)(270,43)(275,45)
\put(285,77){\vector(-1,1){7}}
\put(35,120){\vector(2,3){43}}
\put(145,120){\vector(-2,3){43}}
\put(155,120){\vector(2,3){43}}
\put(265,120){\vector(-2,3){43}}
\put(55,173){\shortstack{$CA$}}
\put(47,163){\shortstack{$CB$}}
\put(39,153){\shortstack{$CC$}}
\put(33,143){\shortstack{$BA$}}
\put(28,133){\shortstack{$AA$}}
\put(22,123){\shortstack{$AC$}}
\put(232,173){\shortstack{$BB$}}
\put(240,163){\shortstack{$BA$}}
\put(248,153){\shortstack{$BC$}}
\put(253,143){\shortstack{$CB$}}
\put(260,133){\shortstack{$AB$}}
\put(266,123){\shortstack{$AC$}}
\put(85,80){\shortstack{$BB$}}
\put(70,112){\shortstack{$CC,\ AA$}}
\put(195,80){\shortstack{$AB,\ BC$}}
\put(205,112){\shortstack{$CA$}}
\put(0,40){\shortstack{$AB,$}}
\put(25,32){\shortstack{$BC$}}
\put(240,40){\shortstack{$AA,$}}
\put(265,32){\shortstack{$CC$}}
\put(85,58){\shortstack{$BA$}}
\put(85,48){\shortstack{$CB$}}
\put(156,53){\shortstack{$AC$}}
\put(104,146){\shortstack{$CA$}}
\put(179,146){\shortstack{$BB$}}
\put(81,208){\shortstack{{\tiny STOP}}}
\put(80,200){\shortstack{{\tiny Alt. 1}}}
\put(201,208){\shortstack{{\tiny STOP}}}
\put(200,200){\shortstack{{\tiny Alt. 2}}}
\put(26,47){\shortstack{$2^n$}}
\put(106,37){\shortstack{$2^{n-1}$}}
\put(184,37){\shortstack{$2^n$}}
\put(93,91){\shortstack{$2^n$}}
\put(74,100){\shortstack{$2^{n-1}$}}
\put(213,91){\shortstack{$2^{n-1}$}}
\put(194,100){\shortstack{$2^n$}}
\put(266,47){\shortstack{$2^n$}}
\put(67,163){\shortstack{\tiny$f_C(x)+f_A(y)$}}
\put(187,163){\shortstack{\tiny$f_B(x)+f_B(y)$}}
\put(230,300){\circle{30}}
\qbezier(220,282)(200,255)(225,245)
\qbezier(240,282)(260,255)(235,245)
\qbezier(225,245)(230,243)(235,245)
\put(245,277){\vector(-1,1){7}}
\put(219,303){\shortstack{{\tiny START}}}
\put(221,295){\shortstack{\tiny $d=0$}}
\put(253,268){\shortstack{$TT$}}
\put(253,258){\shortstack{$LL$}}
\put(253,248){\shortstack{$RR$}}
\put(20,330){\shortstack{$A=T$}}
\put(20,320){\shortstack{$B=L$}}
\put(20,310){\shortstack{$C=R$}}
\put(75,330){\shortstack{$A=R$}}
\put(75,320){\shortstack{$B=L$}}
\put(75,310){\shortstack{$C=T$}}
\put(130,330){\shortstack{$A=L$}}
\put(130,320){\shortstack{$B=T$}}
\put(130,310){\shortstack{$C=R$}}
\put(20,280){\shortstack{$A=R$}}
\put(20,270){\shortstack{$B=T$}}
\put(20,260){\shortstack{$C=L$}}
\put(75,280){\shortstack{$A=T$}}
\put(75,270){\shortstack{$B=R$}}
\put(75,260){\shortstack{$C=L$}}
\put(130,280){\shortstack{$A=L$}}
\put(130,270){\shortstack{$B=R$}}
\put(130,260){\shortstack{$C=T$}}
\put(36,323){\circle{50}}
\put(91,323){\circle{50}}
\put(146,323){\circle{50}}
\put(36,273){\circle{50}}
\put(91,273){\circle{50}}
\put(146,273){\circle{50}}
\put(212,307){\vector(-4,1){43}}
\put(212,293){\vector(-4,-1){43}}
\qbezier(215,310)(150,375)(102,346)
\put(107,348){\vector(-2,-1){7}}
\qbezier(215,290)(150,225)(102,251)
\put(107,249){\vector(-2,1){7}}
\qbezier(218,286)(170,233)(150,233)
\put(150,233){\line(-1,0){60}}
\qbezier(90,233)(50,233)(43,251)
\put(47,245){\vector(-2,3){4}}
\qbezier(218,314)(170,367)(150,367)
\put(150,367){\line(-1,0){60}}
\qbezier(90,367)(50,367)(43,349)
\put(47,355){\vector(-2,-3){4}}
\put(128,309){\line(-1,-1){10}}
\put(128,309){\vector(-1,-1){7}}
\put(128,289){\line(-1,1){10}}
\put(128,289){\vector(-1,1){7}}
\put(73,309){\line(-1,-1){10}}
\put(73,309){\vector(-1,-1){7}}
\put(73,289){\line(-1,1){10}}
\put(73,289){\vector(-1,1){7}}
\put(18,309){\line(-1,-1){10}}
\put(18,309){\vector(-1,-1){7}}
\put(18,289){\line(-1,1){10}}
\put(18,289){\vector(-1,1){7}}
\put(118,299){\line(-1,0){112}}
\qbezier(6,299)(0,299)(0,293)
\put(0,293){\line(0,-1){187}}
\put(0,293){\vector(0,-1){70}}
\qbezier(6,100)(0,100)(0,106)
\put(6,100){\vector(1,0){5}}
\put(123,223){\shortstack{$RL$}}
\put(123,368){\shortstack{$TR$}}
\put(113,247){\shortstack{$TL$}}
\put(113,342){\shortstack{$RT$}}
\put(174,288){\shortstack{$LT$}}
\put(174,304){\shortstack{$LR$}}
\put(4,220){\shortstack{$d=1$}}
\put(-52,210){\shortstack{\tiny deterministic}}
\put(-52,201){\shortstack{\tiny transition}}
\put(-52,192){\shortstack{\tiny (no input read)}}
\put(18,94){\shortstack{\tiny (Alt. 1)}}
\put(141,94){\shortstack{\tiny (draw)}}
\put(258,94){\shortstack{\tiny (Alt. 2)}}
\end{picture}

\vspace{-20.0 pt}
\noindent
Figure 4: The finite state machine: computing $d(x,y)$, the general
case. Follow instructions as for Figure 3. For the deterministic
transition, do not read input or decrease $n$.

\subsection*{4. Translating between the Hanoi graph and $SG_n$}

We now define the graph of states in the $n$-disc Tower of Hanoi
puzzle, and show that it is isomorphic to $SG_n$. The isomorphism may
be computed by reading sequentialy the locations of the discs,
starting with the largest one (which corresponds to the most
significant digit in the Sierpi\'nski gasket labeling), and following a
diagram of permutations translating the labels of the three pegs into
the symbols $T,L,R$ (another finite-state machine!). Together with
the results of the previous section, this will give an effective means
of computing the length of the shortest path between any two vertices
in the Hanoi graph, and of deciding whether the largest disc will be
moved once or twice in a shortest path. After that, we describe
briefly an algorithm for actually constructing the shortest path,
based on the algorithm for getting to a perfect state.

Label the three pegs in the Tower of Hanoi with the symbols
$0,1,2$. Since in a legal state, on each of the pegs the discs are
arranged in increasing size from top to bottom, a state is described
uniquely by specifying, for any disc, the label of its peg. Thus, we
define $H_n$, the $n$-th \emph{Hanoi graph}, to be the graph whose
vertex set is the set $V(H_n) = \{0,1,2\}^n$ (with the coordinates of
the vectors specifying, from left to right, the labels of the pegs of
the largest disc, second-largest disc, etc.), and where edges between
states correspond to permissible moves. Figure 5 shows the graph
$H_4$.

\hspace{-50.0 pt}
\includegraphics{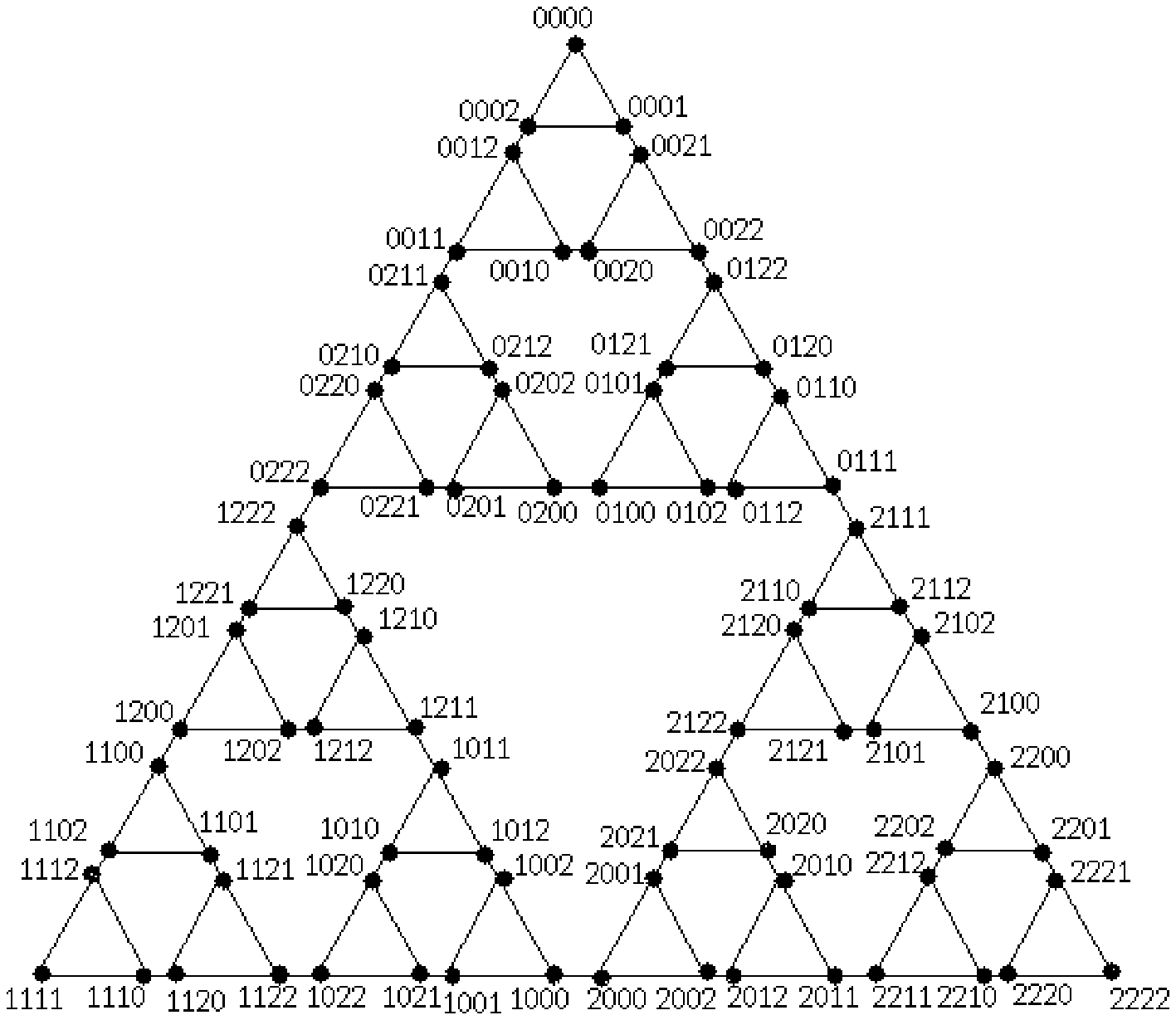}

\vspace{-20.0 pt}
\begin{center} Figure 5: The graph $H_4$ \end{center}

\newpage
The isomorphism between $H_n$ and $SG_n$ is now described by the
following theorem:

\paragraph{Theorem 2.} $H_n$ and $SG_n$ are isomorphic graphs. The
finite-state machine shown below translates a Hanoi state $s \in
\{0,1,2\}^n$ into a Sierpi\'nski gasket labeling $z \in \{T,L,R\}^n$, by
reading the digits from left to right and outputting the symbols
$T,L,R$ at each step according to the identifications in its internal
state, then changing the internal state according to the input.

\begin{center}
\begin{picture}(300,0)(0,135)
\put(22,70){\shortstack{$0=T$}}
\put(22,60){\shortstack{$1=L$}}
\put(22,50){\shortstack{$2=R$}}
\put(36,63){\circle{50}}
\put(22,86){\shortstack{\tiny START}}
\put(96,3){\circle{50}}
\put(82,10){\shortstack{$0=L$}}
\put(82,0){\shortstack{$1=T$}}
\put(82,-10){\shortstack{$2=R$}}
\put(54,45){\vector(1,-1){25}}
\put(60,39){\vector(-1,1){7}}
\put(96,123){\circle{50}}
\put(82,130){\shortstack{$0=T$}}
\put(82,120){\shortstack{$1=R$}}
\put(82,110){\shortstack{$2=L$}}
\put(54,81){\vector(1,1){25}}
\put(60,87){\vector(-1,-1){7}}
\put(170,3){\circle{50}}
\put(156,130){\shortstack{$0=L$}}
\put(156,120){\shortstack{$1=R$}}
\put(156,110){\shortstack{$2=T$}}
\put(170,123){\circle{50}}
\put(156,10){\shortstack{$0=R$}}
\put(156,0){\shortstack{$1=T$}}
\put(156,-10){\shortstack{$2=L$}}
\put(230,63){\circle{50}}
\put(216,70){\shortstack{$0=R$}}
\put(216,60){\shortstack{$1=L$}}
\put(216,50){\shortstack{$2=T$}}
\put(212,45){\vector(-1,-1){25}}
\put(206,39){\vector(1,1){7}}
\put(212,81){\vector(-1,1){25}}
\put(206,87){\vector(1,-1){7}}
\put(119,123){\vector(1,0){28}}
\put(125,123){\vector(-1,0){5}}
\put(119,3){\vector(1,0){28}}
\put(125,3){\vector(-1,0){5}}
\put(60,63){\vector(1,0){146}}
\put(67,63){\vector(-1,0){7}}
\put(106,26){\vector(2,3){50}}
\put(110,32){\vector(-2,-3){6}}
\put(106,100){\vector(2,-3){50}}
\put(110,94){\vector(-2,3){6}}
\put(60,95){\shortstack{$0$}}
\put(200,95){\shortstack{$2$}}
\put(60,24){\shortstack{$2$}}
\put(200,24){\shortstack{$0$}}
\put(90,65){\shortstack{$1$}}
\put(105,84){\shortstack{$2$}}
\put(155,84){\shortstack{$0$}}
\put(130,125){\shortstack{$1$}}
\put(130,-6){\shortstack{$1$}}
\put(0,-36){\shortstack{Figure 6: Computing the isomorphism between
$H_n$ and $SG_n$}}
\end{picture}
\end{center}

\vspace{160.0 pt}

\paragraph{Proof.} This is Lemma 2 in \cite{hinzschief}. There it was claimed
simply that $H_n$ and $SG_n$ are isomorphic, but the proof, which is
by induction, actually describes how to compute the isomorphism, and
this is easily seen to be equivalent to our finite-state machine
formulation.  \qed

\paragraph{Summary.} By running the machines of Figures 4 and 6 in
parallel, we now have an algorithm for computing $d(x,y)$ for two
arbitrary states in the Hanoi graph, and for solving the decision
problem for the largest disc, i.e. to decide whether the largest disc
which it is necessary to move will move once or twice. As we will show
in the next section, when $x$ and $y$ are randomly chosen states, the
expected stopping time of the machine is $63/38$. (This random
variable even has an exponential tail distribution, so with very high
probability only a small number of discs will need to be read to solve
the decision problem.) Having solved the decision problem, the
shortest path may now be computed in a straightforward manner, as
described in \cite{hinz1}, using the algorithm for getting to a
perfect state (use the algorithm described in \cite{hinz1}, or the
algorithm for the Sierpi\'nski gasket described in section 2 together
with the machine of Figure 6 -- which incidentally leads to an
algorithm for getting to perfect states which we have not found in the
literature).

\subsection*{5. The case of random inputs}

\subsubsection*{5.1. How many discs must be read to solve the decision
problem?}

In this section, we calculate the average number of discs that must be
read in order to decide whether in a shortest path the largest disc
will be moved once or twice. Let $x=a_{n-1} a_{n-2}...a_0 \in V(H_n)$,
$y=b_{n-1} b_{n-2} ... b_0 \in V(H_n)$. Assume that we have already
discarded the largest discs which for $x$ and $y$ were on the same
peg, so that $a_{n-1} \neq b_{n-1}$. The algorithm for solving the
decision problem then tells us to run the machines of Figures 4 and 6
until they reach a terminal state (or we run out of input). Since we
have already initialized by discarding irrelevant discs, we will
really be using the machine of Figure 2 (keeping track of the correct
identification of the symbols $L,T,R$ with the pegs $0,1,2$). Since we
are dealing with random inputs, what we are really interested in is
the absorption time of the Markov chain whose transition matrix is
$$
\begin{array}{l} 1 \\ 2 \\ 3 \\ 4 \\ 5 \end{array}
\left( \begin{array}{lllll} 2/9 & 1/9 & 0 & 2/3 & 0 \\
  2/9 & 1/3 & 2/9 & 1/9 & 1/9 \\
  0 & 1/9 & 2/9 & 0 & 2/3 \\
  0 & 0 & 0 & 1 & 0 \\
  0 & 0 & 0 & 0 & 1 \end{array} \right)
$$
into the terminal states 4 and 5. We may identify these two states to
get the simpler matrix
$$
\begin{array}{l} 1 \\ 2 \\ 3 \\ (45) \end{array}
\left( \begin{array}{llll} 2/9 & 1/9 & 0 & 2/3 \\
  2/9 & 1/3 & 2/9 & 2/9 \\
  0 & 1/9 & 2/9 & 2/3 \\
  0 & 0 & 0 & 1 \end{array} \right)
$$
For $i=1,2,3$, denote by $t_i$ the expected time to get to state $(45)$,
\emph{starting from state $i$}. Then clearly we have the equations
$$ t_1 = 1 + \frac{2}{9}t_1 + \frac{1}{9}t_2 $$
$$ t_2 = 1 + \frac{2}{9} t_1 + \frac{1}{3}t_2 + \frac{2}{9}t_3 $$
$$ t_3 = 1 + \frac{1}{9}t_2 + \frac{2}{9}t_3 $$
It may easily be verified that the solution to this system of
equations is
$$ t_1 = 63/38, \quad t_2 = 99/38,\quad t_3 = 63/38 $$
The value $t_1 = 63/38$ is our expected stopping time, since $i=1$
corresponds to the initial state. Note that this value is the limit as
$n\to\infty$ of the average number of discs that must be read; in
reality, for finite $n$ the value will be slightly smaller since after
$n$ steps we run out of input and the machine terminates even if it
has not reached a terminal state. To summarize:

\paragraph{Theorem 3.} The decision problem for shortest paths can be
solved in average time $O(1)$. Specifically, the average number of
disc pairs that our algorithm must read, once identical discs have been
discarded, is bounded from above by, and converges as $n\to\infty$ to,
$63/38$.

\subsubsection*{5.2. The average distance between points on the
Sierpi\'nski gasket}

Hinz and Schief \cite{hinzschief} computed the average length
$466/885$ of a shortest path between two random points on the
\emph{infinite} Sierpi\'nski gasket of unit side. An equivalent result
of Hinz \cite{hinz3} and of Chan \cite{chan}, in terms of the
Tower of Hanoi, is that the average number of moves in a shortest path
between two random states in the $n$-disc Tower of Hanoi, is
asymptotically $(466/885)\cdot 2^n$ as $n\to\infty$.

Without going into too much detail, we show that it is
possible to obtain the value of $466/885$ just by looking at the
finite-state machine of Figure 4. Since we are dealing with the
infinite gasket, we start with $n=0$ and, as before, decrease the
value of $n$ after each step, so that $n$ will go into the negative
integers. Let $d_1, d_2, d_3, d_4$ be the expected accumulated values
of the variable $d$ if one starts the machine, with initial values
$n=0, d=0$, at either of the four non-terminal states, in order of
their distance from the state {\tiny START} (so $d_1$ is the total
distance; $d_2$ is the distance after discarding identical
most-significant digits of $x$ and $y$, etc.). Then we have the
equations
$$ d_1 = \frac{1}{3}\cdot \frac{1}{2}d_1 + \frac{2}{3}\cdot
\frac{1}{2}d_2 $$
$$ d_2 = \frac{2}{9}\cdot \left(1+\frac{1}{2}d_2\right) +
\frac{1}{9}\cdot \left(1+\frac{1}{2}d_3\right) +
\frac{2}{3}\cdot \left(\frac{1}{2}+\frac{2}{3}\right) $$
$$ d_3 = \frac{2}{9}\cdot\left(\frac{1}{2}+\frac{1}{2}d_2\right) +
\frac{2}{9}\cdot \left(\frac{1}{2}+\frac{1}{2}d_3\right) +
\frac{1}{9}\cdot \left(1+\frac{1}{2}d_3\right) + $$ $$
\frac{2}{9}\cdot\left(\frac{1}{2}+\frac{1}{2}d_4\right) + 
\frac{2}{9}\cdot \frac{2}{3} $$
$$ d_4 = \frac{1}{9}\cdot \left(1+\frac{1}{2}d_3\right) +
\frac{2}{9}\cdot \left(1+\frac{1}{2}d_4\right) + \frac{2}{3}\cdot
\left(\frac{1}{2}+\frac{2}{3}\right) $$ 

The value (1/2+2/3) in the second and fourth equations is the expected
value of $f_C(x)+f_A(y)$ (respectively $f_B(x)+f_B(y)$), given that
the first pair of inputs is one of the six values $AC, AA, BA, CC, CB,
CA$ (respectively $BB, BA, BC, CB, AB, AC$).

Again, it may verified that the solution to this system of equations
is
$$ d_1 = \frac{466}{885}, \quad d_2=\frac{233}{177},\quad 
   d_3= \frac{188}{177},\quad d_4=\frac{233}{177} $$
which gives our claimed value for $d_1$.

\end{document}